\documentclass[a4paper,11pt]{article}
\usepackage{amssymb} \usepackage{graphicx} \usepackage{color}
\usepackage[utf8]{inputenc} %set che funziona su Ubuntu al posto di latin1
\usepackage[T1]{fontenc} \usepackage[english]{babel}
 \usepackage{amsmath,amsthm} \usepackage{makeidx}
\newcommand{\tr}{\textcolor{red}} \newcommand{\mb}{\mbox{}} \newcommand{\mbs}{\mbox{ }}
\newcommand{\bn}{b_n := \frac{a_n}{a_{n-1}} - \frac{a_{n+1}}{a_n}} \newcommand{\bnr}{\frac{b_{n-1}}{b_n}}
\newcommand{\cdotm}{\mbs \cdot \mbs} \newcommand{\minusm}{\mbs - \mbs}

\newcommand{\suppr}{\hfuzz=100pt \vfuzz=100pt} \suppr
\usepackage[numbers,sort&compress]{natbib}

\begin{document}

\newtheorem{theorem}{Teorema}[section]
\newtheorem{definition}[theorem]{Definition}

\title{On certain ratios regarding integer numbers which are both triangulars and squares}
\author{Fabio Roman\footnote{University of Torino, Dipartimento di Matematica, Via Carlo Alberto 10, 10123 Torino. Tel: +39 011 6702827; Fax: +39 011 6702878; e-mail: fabio.roman@unito.it}}
\date{}

\maketitle

\begin{abstract}
We investigate integer numbers which possess at the same time the properties to be triangulars and squares, that are, numbers $a$ for which do exist integers $m$ and $n$ such that $ a = n^2 = \frac{m \cdot (m+1)}{2} $. In particular, we are interested about ratios between successive numbers of that kind. While the limit of the ratio for increasing $a$ is already known in literature, to the best of our knowledge the limit of the ratio of differences of successive ratios, again for increasing $a$, is a new investigation. We give a result for the latter limit, showing that it coincides with the former one, and we formulate a conjecture about related limits.
\end{abstract}

\section{Preliminaries}
We recall some basic definitions from elementary number theory.
\begin{definition}
A non-negative integer is said to be \textbf{triangular} if it can be the number of objects in a set able to form a triangle, right or equilateral.

%\ section{Preliminari}
%Prima di cominciare la trattazione vera e propria, vale la pena ricordare le definizioni di base di alcuni oggetti matematici che saranno considerati:
%\begin{definition} 
%Un numero si definisce \textbf{triangolare} se può essere la quantità di un insieme di oggetti disposti in modo da formare un triangolo, rettangolo o equilatero. \\

$$ \begin{tabular}{ccccc}
$\bullet$ & & & & \\
$\bullet$ & $\bullet$ & & & \\
$\bullet$ & $\bullet$ & $\bullet$ & & \\
$\bullet$ & $\bullet$ & $\bullet$ & $\bullet$ & \\
$\bullet$ & $\bullet$ & $\bullet$ & $\bullet$ & $\bullet$
\end{tabular} \\ \\ $$

\noindent A triangular number has the form $$ T_n := \frac{n \cdot (n+1)}{2} = \binom{n+1}{2} $$ where $n$ is a natural number.

%\noindent Un numero triangolare ha la forma: $$ T_n := \frac{n \cdot (n+1)}{2} = \binom{n+1}{2} $$ dove $n$ è un numero naturale. 

\end{definition}

\begin{definition}

Similarly, a non-negative integer is said to be \textbf{square} if it can be the number of objects in a set able to form a square.
 
%Un numero si definisce \textbf{quadrato}, in maniera analoga alla definizione precedente, se può essere la quantità di un insieme di oggetti disposti in modo da formare un quadrato. \\

$$ \begin{tabular}{ccccc}
$\bullet$ & $\bullet$ & $\bullet$ & $\bullet$ & $\bullet$ \\
$\bullet$ & $\bullet$ & $\bullet$ & $\bullet$ & $\bullet$ \\
$\bullet$ & $\bullet$ & $\bullet$ & $\bullet$ & $\bullet$ \\
$\bullet$ & $\bullet$ & $\bullet$ & $\bullet$ & $\bullet$ \\
$\bullet$ & $\bullet$ & $\bullet$ & $\bullet$ & $\bullet$
\end{tabular} \\ \\ $$

%Resta inteso che i numeri quadrati sono tutti e soli quelli della forma $ n \cdot n = n^2 $, con $n$ numero naturale.

\noindent It is straightforward to say that square numbers have the form $ n \cdot n = n^2 $, where $n$ is a natural number.

\end{definition}

%\noindent Analogamente alla definizione di numero triangolare e di numero quadrato, è possibile definire generalmente un \emph{numero poligonale}, se può essere la quantità di un insieme di oggetti disposti in modo da formare un poligono regolare di un numero determinato di lati. Nella fattispecie:
%\begin{definition} Un numero si definisce \textbf{$m$-gonale} se può essere la quantità di un insieme di oggetti disposti in modo da formare un poligono regolare di $m$ lati. Sono tutti e soli i numeri della forma:\ label{mgonale}
%$$ P_{m,n} = \frac{n[(m-2)n-(m-4)]}{2} $$ \end{definition}

\noindent We can also define a generic \emph{polygonal number} as an integer that can be the number of objects in a set able to form a regular polygon having a certain number of sides.
\begin{definition} A number is said to be \textbf{$m$-gonal} if it can be the number of objects in a set able to form a regular $m$-gon. The $n$-th $m$-gonal number has the form:
\label{mgonale} $$ P_{m,n} = \frac{n[(m-2)n-(m-4)]}{2} $$ \end{definition}

%\newpage

%\ section{I calcoli di base} \ label{calcolibase}
%Facendo riferimento alle definizioni di numero quadrato e triangolare, appena riportate, e impostando l'uguaglianza, si ottiene:
%$$ n^2 = \frac{m \cdot (m+1)}{2} $$
%Una successione di passaggi algebrici porta a:
%\begin{eqnarray*}
%n^2 & = & \frac{m \cdot (m+1)}{2} \\
%n^2 & = & \frac{m^2+m}{2} \\
%2n^2 - (m^2+m) & = & 0 \\
%2n^2 - \left( m^2+m+\frac{1}{4} \right) & = & -\frac{1}{4} \\
%2n^2 - \left( m+\frac{1}{2} \right)^2 & = & -\frac{1}{4} \\
%8n^2 - (2m+1)^2 & = & -1 \\
%t^2 - 8n^2 & = & 1
%\end{eqnarray*}
%dove si è posto infine $ t := 2m+1 $. \\
%Le soluzioni sono dunque tutte e sole le coppie $(t,n)$ dell'equazione di Pell cui sopra con il vincolo $t$ dispari.
%\\ È possibile, volendo, porre $ s:= 2n $ e far diventare l'equazione $ t^2 - 2s^2 = 1$, riconducedosi alla più classica delle equazioni di Pell, con tuttavia il vincolo supplementare della parità di $s$. \\
%Come ben noto dalla teoria, esistono infinite coppie di numeri interi positivi che soddisfano l'equazione, ed esistono formule ricorsive e formule chiuse che permettono di determinare queste coppie al variare di un parametro, che vedremo di esplicitare in seguito. \ label{pell}

\section{Basic computations} 
\label{calcolibase}
With reference to square and triangular numbers' definitions, by imposing equality, we obtain:
$$ n^2 = \frac{m \cdot (m+1)}{2} $$
that can be algebrically transformed to:
\begin{eqnarray*}
n^2 & = & \frac{m \cdot (m+1)}{2} \\
n^2 & = & \frac{m^2+m}{2} \\
2n^2 - (m^2+m) & = & 0 \\
2n^2 - \left( m^2+m+\frac{1}{4} \right) & = & -\frac{1}{4} \\
2n^2 - \left( m+\frac{1}{2} \right)^2 & = & -\frac{1}{4} \\
8n^2 - (2m+1)^2 & = & -1 \\
t^2 - 8n^2 & = & 1
\end{eqnarray*}
by setting at the end $ t := 2m+1 $. \\
This allows us to say that $(t,n)$ should solve a Pell equation, assuming that $t$ is odd; we can also set $ s := 2n $, in order to write $ t^2 - 2s^2 = 1 $, that is the more classical Pell equation, in which we need $s$ even. \label{pell} 

%\newpage % bocconi/cesenatico, tesi dottorato, mamma con lavoro

%\ section{La trattazione empirica}
%Vediamo come si può trattare questo problema empiricamente con l'utilizzo di un semplice foglio di calcolo (spreadsheet), strumento di uso molto comune presente pressoché su ogni personal computer. \\
%L'idea è quella di creare una tabella di questo tipo. Nella prima colonna andiamo ad inserire i numeri interi positivi da $1$ fino al numero desiderato; nella seconda, facciamo calcolare l'$n$-esimo numero triangolare, dove $n$ è il numero della riga considerata; nella terza, viene calcolata la radice quadrata di tale numero. Nella quarta riga tale radice viene arrotondata per difetto (ovvero ne viene presa la parte intera), mentre nella quinta riga si fa la differenza tra la terza e la quarta (parte decimale o mantissa). Esplicitando: \\

\section{A numerical approach}
We treat here the problem empirically, by using a spreadsheet. \\
The idea is to create a table, where in the first column is listed a certain number of positive integers, in the second one the respective triangular number, in the third one its square root; in the fourth column we chop the square root at the lower integer, while in the fifth and last column, we do the difference between the third one and the fourth one, obtaining its decimal part. The first ten rows of the table give:

$$ \begin{tabular}{|c|c|c|c|c|}
\hline \textbf{integers} & \textbf{triangulars} & \textbf{roots} & \textbf{integer parts} & \textbf{decimal parts} \\
\hline 1 & \tr{1} & 1.0000 & 1.0000 & \tr{0.0000} \\
\hline 2 & 3 & 1.7321 & 1.0000 & 0.7321 \\
\hline 3 & 6 & 2.4495 & 2.0000 & 0.4495 \\
\hline 4 & 10 & 3.1623 & 3.0000 & 0.1623 \\
\hline 5 & 15 & 3.8730 & 3.0000 & 0.8730 \\
\hline 6 & 21 & 4.5826 & 4.0000 & 0.5826 \\
\hline 7 & 28 & 5.2915 & 5.0000 & 0.2915 \\
\hline 8 & \tr{36} & 6.0000 & 6.0000 & \tr{0.0000} \\
\hline 9 & 45 & 6.7082 & 6.0000 & 0.7082 \\
\hline 10 & 55 & 7.4162 & 7.0000 & 0.4162 \\ \hline
\end{tabular} $$

%e via di seguito. \\
%Risulta pressoché immediato affermare che il numero triangolare in questione è anche un quadrato se e soltanto se il risultato in quinta colonna è pari a $0$ (per numeri particolarmente grandi potrebbero iniziare ad insorgere problemi dovuti all'aritmetica finita di macchina, ma questioni di questo {tipo} non costituiscono una problematica in questo caso, considerati gli ordini di grandezza trattati). \\
%Proseguendo con la tabella fino a $65534$ (la versione dello spreadsheet utilizzata allora disponeva di $2^{16}-1 = 65535$ righe, e la prima riga veniva utilizzata come intestazione delle colonne), è possibile trovare per via diretta altri di questi numeri:

\noindent and so on. \\
It is straightforward to say that the considered triangular number is also a squadre if and only if, for its row, the fifth column is 0. Exceptions can arise due to finite arithmetic errors, but it's not the case at least for the moment, because all numbers appearing are not too large to generate machine-caused loss of precision. \\
By extending with the table until $65534$ (if we use a spreadsheet with $2^{16}-1 = 65535$ rows, using the first for naming columns), we can directly found some of these numbers:

$$ \begin{tabular}{|c|c|c|c|c|}
\hline \textbf{integers} & \textbf{triangulars} & \textbf{roots} & \textbf{integer parts} & \textbf{decimal parts} \\
\hline 49 & \tr{1225} & 35.0000 & 35.0000 & \tr{0.0000} \\
\hline 288 & \tr{41616} & 204.0000 & 204.0000 & \tr{0.0000} \\ 
\hline 1681 & \tr{1 \mb 413 \mb 721} & 1189.0000 & 1189.0000 & \tr{0.0000} \\
\hline 9800 & \tr{48 \mb 024 \mb 900} & 6930.0000 & 6930.0000 & \tr{0.0000} \\
\hline 57121 & \tr{1 \mb 631 \mb 432 \mb 881} & 40391.0000 & 40391.0000 & \tr{0.0000} \\ \hline
\end{tabular} $$

%\newpage

%\noindent Una proprietà che balza all'occhio è come il rapporto tra due numeri successivi ad essere sia triangolari che quadrati sia \textit{all'incirca lo stesso} nei vari casi che si succedono. Analizzando da vicino questa cosa:

\noindent A property that can be observed is that the ratio between two successive numbers, both triangular and square, seems to be the same, case after case. By explicit computation:

$$ \begin{tabular}{|c|c|}
\hline \textbf{$a_n$} & \textbf{$a_{n+1}/a_n$} \\
\hline 1 & 36.00000 \\
\hline 36 & 34.02778 \\
\hline 1225 & 33.97224 \\
\hline 41616 & 33.97061 \\
\hline 1 \mb 413 \mb 721 & 33.97056420609 \\
\hline 48 \mb 024 \mb 900 & \\ \hline
\end{tabular} $$

%si vede che il rapporto tende piuttosto in fretta ad una data quantità. Ma si può dire di più, ed è il motivo per cui alla fine abbiamo tenuto 11 cifre decimali invece delle 5 tenute in precedenza. Il rapporto tra le differenze dei rapporti tende anch'esso alla stessa quantità. Infatti:

\noindent we can see how this ratio seems to rapidly converge to a fixed value. But we can say more of that, and this is why we kept 11 digits instead of 5 in the last step: also the ratio between differences of subsequent ratios converge to the same value. In fact:

$$ \begin{tabular}{|c|c|c|c|}
\hline \textbf{$a_n$} & \textbf{$\frac{a_{n+1}}{a_n}$} & \textbf{$\bn$} & \textbf{$\bnr$} \\
\hline 1 & 36.00000 & & \\
\hline 36 & 34.02778 & 1.97222 & \\
\hline 1225 & 33.97224 & 0.05553 & 35.51450 \\
\hline 41616 & 33.97061 & 0.00163 & 34.01430 \\
\hline 1 \mb 413 \mb 721 & 33.97056420609 & 0.0000480584221 & 33.97185 \\
\hline 48 \mb 024 \mb 900 & & & \\ \hline
\end{tabular} $$

%Allora, per andare avanti con la tabella, proviamo a procedere nel seguente modo:
%\begin{itemize}
%\item dividiamo l'ultimo valore della successione $b_n$ riscontrato, pari a $0.0000480584221$, per l'analogo della successione di rapporti che nella tabella sta alla sua destra, ovvero $33.97185$;
%\item otteniamo la quantità $0.0000014146543$, che andiamo a togliere dall'ultimo valore disponibile del rapporto in seconda colonna, ovvero $33.97056420609$;
%\item otteniamo la quantità $33.97056279144$, che andiamo a moltiplicare per l'ultimo numero disponibile, $ 48 \mbs 024 \mbs 900 $; se quanto visto prima continua ad avere luogo, si dovrebbe ottenere un numero (reale, non necessariamente intero) che approssima piuttosto bene un nuovo numero sia triangolare che quadrato.
%\end{itemize}
%Il risultato è sorprendente:
%$$ 48 \mbs 024 \mbs 900 \cdot 33.97056279144 = 1 \mbs 631 \mbs 432 \mbs 881.00263 $$
%e $ 1 \mbs 631 \mbs 432 \mbs 881 $ è effettivamente un numero sia triangolare che quadrato, come abbiamo potuto verificare direttamente con una delle tabelle precedenti. \\
%Utilizzando i valori appena ricavati, per ognuna delle quantità in questione, possiamo aggiungere una riga alla tabella:

So, we can create new rows in the table, by following these steps:
\begin{itemize}
\item we divide the latest value $b_n$ obtained, i.e. $0.0000480584221$, for the analogous from the sequence of ratios which in the table lies on its right, so $33.97185$;
\item we obtain $0.0000014146543$, and we subtract it from the latest value available in the second column, specifically $33.97056420609$;
\item we obtain $33.97056279144$, and we multiply it for the latest number written, $ 48 \mbs 024 \mbs 900 $; if what we conjectured is correct, we should obtain a number (real, not necessarily integer) well-approximating a new both triangular and square number.
\end{itemize}
In fact, we compute:
$$ 48 \mbs 024 \mbs 900 \cdot 33.97056279144 = 1 \mbs 631 \mbs 432 \mbs 881.00263 $$
and $ 1 \mbs 631 \mbs 432 \mbs 881 $ is both triangular and square. That allows us to add a line into the table:

$$ \begin{tabular}{|c|c|c|c|}
\hline \textbf{$a_n$} & \textbf{$\frac{a_{n+1}}{a_n}$} & \textbf{$\bn$} & \textbf{$\bnr$} \\
\hline 1 & 36.00000 & & \\
\hline 36 & 34.02778 & 1.97222 & \\
\hline 1225 & 33.97224 & 0.05553 & 35.51450 \\
\hline 41616 & 33.97061 & 0.00163 & 34.01430 \\
\hline 1 \mb 413 \mb 721 & 33.97056420609 & 0.0000480584221 & 33.97185 \\
\hline 48 \mb 024 \mb 900 & 33.97056279144 & 0.0000014146543 & 33.97185 \\
\hline 1 \mb 631 \mb 432 \mb 881.00263 & & & \\ \hline
\end{tabular} $$

%Quello che abbiamo utilizzato in questo caso è un procedimento \textit{all'indietro}. Ovvero, osservando che la quantità dell'ultima colonna tende a stabilizzarsi, abbiamo stimato, accettando un certo errore, che si stabilizzasse esattamente all'ultimo valore noto, calcolando poi di conseguenza tutti i valori della stessa riga e delle colonne precedenti (tranne la prima colonna che è sfasata di una riga e della quale si calcola quindi la riga successiva, ma è una scelta meramente tipografica). \\
%Possiamo però ora rettificare la tabella, tenuto conto del fatto che conosciamo il valore esatto del nuovo numero sia triangolare che quadrato, e quindi possiamo passare dal procedimento all'indietro a quello in avanti: dalla conoscenza dei numeri che possiedono tale proprietà, ricavare i rapporti e le differenze.

\noindent We used a \emph{backward} completion: by observing that the value in the last column tends to stabilize, we estimate, by accepting an error margin, that it is constant starting from the considered row, and we complete the row by calculating all values in the previous columns. \\
By taking account of the fact that we know the exact value of the new number both triangular and square, we can rectify the table, moving from backward completion to \emph{forward} completion: if we know the numbers having this property, we can derive ratios and differences. 

$$ \begin{tabular}{|c|c|c|c|}
\hline \textbf{$a_n$} & \textbf{$\frac{a_{n+1}}{a_n}$} & \textbf{$\bn$} & \textbf{$\bnr$} \\
\hline 1 & 36.00000 & & \\
\hline 36 & 34.02778 & 1.97222 & \\
\hline 1225 & 33.97224 & 0.05553 & 35.51450 \\
\hline 41616 & 33.97061 & 0.00163 & 34.01430 \\
\hline 1 \mb 413 \mb 721 & 33.97056420609 & 0.0000480584221 & 33.97185 \\
\hline 48 \mb 024 \mb 900 & 33.97056279139 & 0.0000014147063 & 33.97060 \\
\hline 1 \mb 631 \mb 432 \mb 881 & & & \\ \hline
\end{tabular} $$

%La quantità a cui tende il rapporto di sinistra sembra essere approssimabile, utilizzando cinque cifre decimali, a $33.97056$. Proviamo a prendere tale numero anche come quantità per il rapporto di destra in modo tale da vedere se in questo modo è possibile ricavare dei numeri successivi. Come prima:
%\begin{itemize}
%\item $0.0000014147063$ diviso $33.97056$ uguale $0.0000000416451$;
%\item $33.97056279139$ meno $0.0000000416451$ uguale $33.97056274974$;
%\item $ 1 \mbs 631 \mbs 432 \mbs 881 $ per $33.97056274974$ uguale $ 55 \mbs 420 \mbs 693 \mbs 055.99960 $.
%\end{itemize}
%Si verifica, con un qualunque software di calcolo numerico, che $ 55 \mbs 420 \mbs 693 \mbs 056 $ è:
%\begin{itemize}
%\item il $332928$-esimo triangolare: infatti, detto $c$ il numero, l'equazione \\ $n^2+n-2c=0$ possiede quella come soluzione positiva;
%\item il $235416$-esimo quadrato: l'equazione $n^2=c$ possiede quella come soluzione positiva.
%\end{itemize}
%Possiamo di nuovo prima aggiornare e poi rettificare la tabella:

\noindent It seems we can approximate the limit of the left ratio, using 5 digits, to $33.97056$; if we try to take the same number as the right ratio, hoping to find new numbers, we can proceed:
\begin{itemize}
\item $ 0.0000014147063 / 33.97056 = 0.0000000416451 $;
\item $ 33.97056279139 - 0.0000000416451 = 33.97056274974 $;
\item $ 1 \mbs 631 \mbs 432 \mbs 881 \cdot 33.97056274974 = 55 \mbs 420 \mbs 693 \mbs 055.99960 $
\end{itemize}
Any \texttt{CAS} allows us to consider $ 55 \mbs 420 \mbs 693 \mbs 056 $ as:
\begin{itemize}
\item the $332928$-th triangular number: if $c$ is the number, the algebraic equation $ n^2 + n - 2c = 0 $ has that value of $n$ as positive root;
\item the $235416$-th square number, just by calculating its square root.
\end{itemize}
These data allows us, again, to update and rectify the table:

$$ \begin{tabular}{|c|c|c|c|}
\hline \textbf{$a_n$} & \textbf{$\frac{a_{n+1}}{a_n}$} & \textbf{$\bn$} & \textbf{$\bnr$} \\
\hline 1 & 36.00000 & & \\
\hline 36 & 34.02778 & 1.97222 & \\
\hline 1225 & 33.97224 & 0.05553 & 35.51450 \\
\hline 41616 & 33.97061 & 0.00163 & 34.01430 \\
\hline 1 \mb 413 \mb 721 & 33.97056420609 & 0.0000480584221 & 33.97185 \\
\hline 48 \mb 024 \mb 900 & 33.97056279139 & 0.0000014147063 & 33.97060 \\
\hline 1 \mb 631 \mb 432 \mb 881 & 33.97056274974 & 0.0000000416451 & 33.97056 \\
\hline 55 \mb 420 \mb 693 \mb 055.99960 & & & \\ \hline
\end{tabular} $$

$$ \begin{tabular}{|c|c|c|c|}
\hline \textbf{$a_n$} & \textbf{$\frac{a_{n+1}}{a_n}$} & \textbf{$\bn$} & \textbf{$\bnr$} \\
\hline 1 & 36.00000 & & \\
\hline 36 & 34.02778 & 1.97222 & \\
\hline 1225 & 33.97224 & 0.05553 & 35.51450 \\
\hline 41616 & 33.97061 & 0.00163 & 34.01430 \\
\hline 1 \mb 413 \mb 721 & 33.97056420609 & 0.0000480584221 & 33.97185 \\
\hline 48 \mb 024 \mb 900 & 33.97056279139 & 0.0000014147063 & 33.97060 \\
\hline 1 \mb 631 \mb 432 \mb 881 & 33.97056274974 & 0.0000000416451 & 33.97056 \\
\hline 55 \mb 420 \mb 693 \mb 056 & & & \\ \hline
\end{tabular} $$

%Abbiamo in un certo senso \textit{accelerato} il procedimento, infatti operando come prima avremmo dovuto prendere come valore dell'ultima colonna $33.97060$, ovvero l'ultimo valore disponibile, invece abbiamo preso $33.97056$, nell'ipotesi che, così come i rapporti della seconda colonna tendono a quella quantità, lo stesso facciano quelli della quarta colonna.
%\\ Si noti che effettivamente abbiamo ottenuto un valore così vicino al vero che la rettifica della prima colonna non porta ad alcun cambiamento nelle altre, limitatamente al numero di cifre decimali considerato; si può anche osservare che, come mostra la tabella successiva, tenendo in considerazione sempre lo stesso numero di cifre dopo la virgola, lo stesso risultato si sarebbe ottenuto prendendo $33.97060$ come valore del rapporto di destra, tuttavia aumentando il numero di cifre ci si può aspettare che la precisione sia inferiore:

We have in some sense \emph{fastened} the procedure: in fact, by following what we done before, we would have taken as value in the last column $33.97060$, i.e. the last value available, while we took instead $33.97056$, assuming that ratios in the fourth column converge at the same quantity ratios in the second column do.
\\ We note that we obtain an almost exact value: a rectify in the first column doesn't change anything in the others, with respect to the number of digits considered; we can also observe, as seen in the next table, that by using the same number of digits, we would obtain the same result by taking $33.97060$ as right ratio, while an increment in the number of digits would likely result in a difference, in which the lower precision lies in the choice of that value.

$$ \begin{tabular}{|c|c|c|c|}
\hline \textbf{$a_n$} & \textbf{$\frac{a_{n+1}}{a_n}$} & \textbf{$\bn$} & \textbf{$\bnr$} \\
\hline 1 & 36.00000 & & \\
\hline 36 & 34.02778 & 1.97222 & \\
\hline 1225 & 33.97224 & 0.05553 & 35.51450 \\
\hline 41616 & 33.97061 & 0.00163 & 34.01430 \\
\hline 1 \mb 413 \mb 721 & 33.97056420609 & 0.0000480584221 & 33.97185 \\
\hline 48 \mb 024 \mb 900 & 33.97056279139 & 0.0000014147063 & 33.97060 \\
\hline 1 \mb 631 \mb 432 \mb 881 & 33.97056274974 & 0.0000000416450 & 33.97060 \\
\hline 55 \mb 420 \mb 693 \mb 055.99960 & & & \\ \hline
\end{tabular} $$

By applying again the method, we obtained another couple of numbers: $ 1 \mbs 882 \mbs 672 \mbs 131 \mbs 025 $ and $ 63 \mbs 955 \mbs 431 \mbs 761 \mbs 796 $. \\
We can also note that:
\begin{itemize} \label{rapporti}
\item if we define $c_n := \frac{b_{n-2}}{b_{n-1}} - \frac{b_{n-1}}{b_n}$, even $\frac{c_{n-1}}{c_n}$ tends to the same value; we can conjecture that it happens every time we iterate in this way, that is, if we denote $a_{1,n}:=a_n, \mbs a_{2,n}:=b_n, \mbs a_{3,n}:=c_n$, and we define for every $ i \geq 3 $ a corresponding $a_{i,n} := \frac{a_{i-1,n-2}}{a_{i-1,n-1}} -\frac{a_{i-1,n-1}}{a_{i-1,n}}$, we can say that, again for every $ i \geq 3 $, while $n$ tends to infinity, $\frac{a_{i,n-1}}{a_{i,n}}$ tends to the value.
\item if we use more digits for the ratios, and we assume correct the conjecture, we can consider one of the ratios, call $d$ the difference between a value and the previous one, $q$ the recurring value of about $33.97056$, and say that the subsequent difference will be approximable by $\frac{d}{q}$, the next one by $\frac{d}{q^2}$, and so on. The sum of the difference from there to infinity will be approximable by  $ \frac{d}{q} + \frac{d}{q^2} + \frac{d}{q^3} + \ldots = \frac{d}{q-1} = \frac{d}{32.97056} $, that allows us to obtain a gain in the relative precision of at least $32$ times every single step, and at least $1000$ times every two steps, that corresponds to three digits.
\end{itemize}

%Naturalmente però il metodo comporta anche dei limiti, nella fattispecie:
%\begin{itemize}
%\item la richiesta di una certa precisione computazionale: con $15$ cifre decimali significative, come da standard del tipo \texttt{float}, già a calcolare il numero più grande riportato, intero di $14$ cifre, iniziano a verificarsi imprecisioni; moltiplicando per $q$ tale numero, si ottiene un intero di $16$ cifre che non può quindi essere rappresentato esattamente con il sistema dei reali in doppia precisione a 64 bit.
%\item premesso che, stante in questo modo la trattazione, non c'è il supporto teorico a confermare che si ottengono effettivamente \textit{tutti} numeri che possiedono la proprietà di essere sia triangolari che quadrati, sebbene così pare risultare, non si può verificare che siano i \textit{soli} numeri a possedere tale proprietà, ovvero che, dati due elementi consecutivi della successione (finita o infinita, questo a priori non lo sappiamo), non ci sia un numero compreso tra questi che possiede anch'esso la proprietà.
%\end{itemize}
%Vedremo nel seguito come questo numero reale pari all'incirca a $33.97056$ non sia una quantità che si possa dire essere \textit{nuova} (anzi, trattasi di un numero algebrico e non trascendentale), e come ci sia una solida trattazione matematica di tipo teorico a supportare tutto ciò.

On the other hand, we need a certain machine precision: with $15$ digits, that corresponds to a relative precision of about $2^{-52}$, the standard of the \texttt{double} type, we report a loss of precision in the computation of the biggest number found before, a $14$-digit integer. If we multiply that number for $q$, we obtain a $16$-digit integer, and in general we can't exactly write a $16$-digit integer as a 64-bit real value.

\section{Exact approach with Pell equations}
\label{applicazione}
It is widely known from Pell equations' theory that, for solving:
$$ t^2 - 2s^2 = 1 $$
we start by write $\sqrt{2}$ as a continuous fraction, that is:
$$ \sqrt{2} = 1+\frac{1}{2+\frac{1}{2+\frac{1}{2+\ldots}}} $$
The first convergent is $\frac{3}{2}$, and $(t,s)=(3,2)$ does in fact solve the equation, i.e. $3^2 - 2 \cdot 2^2 = 9 - 8 = 1 $. \\
By the relation $s=2n$, we have $n=1$, and $n^2=1$, that is the first number both triangular and square. \\
Successive integers can be found in a traditional way, involving well-estabilished theory:

\newcommand{\sq}{\sqrt{2}}
$$ \begin{tabular}{|c|c|c|c|c|c|c|}
\hline $i$ & $(3+2\sq)^i$ & $t$ & $s$ & $m$ & $n$ & $n^2$ \\
\hline $1$ & $(3+2\sq)$ & $3$ & $2$ & $1$ & $1$ & $1$ \\
\hline $2$ & $(17+12\sq)$ & $17$ & $12$ & $8$ & $6$ & $36$ \\
\hline $3$ & $(99+70\sq)$ & $99$ & $70$ & $49$ & $35$ & $1225$ \\
\hline $4$ & $(577+408\sq)$ & $577$ & $408$ & $288$ & $204$ & $41616$ \\
\hline $5$ & $(3363+2378\sq)$ & $3363$ & $2378$ & $1681$ & $1189$ & $1413721$ \\
\hline $6$ & $(19601+13860\sq)$ & $19601$ & $13860$ & $9800$ & $6930$ & $48024900$ \\ 
\hline \end{tabular} $$

and again: % e ancora:

$$ \begin{tabular}{|c|c|c|c|c|c|}
\hline $i$ & $t$ & $s$ & $m$ & $n$ & $n^2$ \\
\hline $7$ & $114243$ & $80782$ & $57121$ & $40391$ & $1 \mbs 631 \mbs 432 \mbs 881$ \\
\hline $8$ & $665857$ & $470832$ & $332928$ & $235416$ & $55 \mbs 420 \mbs 693 \mbs 056$ \\
\hline $9$ & $3880899$ & $2744210$ & $1940449$ & $1372105$ & $1 \mbs 882 \mbs 672 \mbs 131 \mbs 025$ \\
\hline $10$ & $22619537$ & $15994428$ & $11309768$ & $7997214$ & $63 \mbs 955 \mbs 431 \mbs 761 \mbs 796$ \\
\hline \end{tabular} $$

%e così via. Osserviamo che abbiamo ricavato tutti e soli i numeri che avevamo ottenuto con il metodo empirico; questa volta però abbiamo alle spalle una solida motivazione teorica (paragrafo \ref{ulteriori}) per affermare che non ci possono essere altri numeri sia triangolari che quadrati compresi tra quelli riportati in due righe consecutive qualsiasi della tabella. \\ \\
%In generale:
%\begin{eqnarray*}
%(t_{i-1} + s_{i-1}\sq)(3+2\sq) & = & (t_i + s_i\sq) \\
%3 t_{i-1} + 2\sq t_{i-1} + 3\sq s_{i-1} + 4 s_{i-1} & = & (t_i + s_i\sq) \\
%3 t_{i-1} + 4 s_{i-1} + (2 t_{i-1} + 3 s_{i-1}\sq) & = & (t_i + s_i\sq) \\
%\end{eqnarray*}
%da cui, ricorsivamente:
%$$ \left\{	\begin{aligned}
%t_i & = & 3 t_{i-1} + 4 s_{i-1} \\
%s_i & = & 2 t_{i-1} + 3 s_{i-1} 
%\end{aligned}	\right.	$$
%%\begin{eqnarray*}
%%t_i & = & 3 t_{i-1} + 4 s_{i-1} \\
%%s_i & = & 2 t_{i-1} + 3 s_{i-1} 
%%\end{eqnarray*}

and so on; we can generalize:
\begin{eqnarray*}
(t_{i-1} + s_{i-1}\sq)(3+2\sq) & = & (t_i + s_i\sq) \\
3 t_{i-1} + 2\sq t_{i-1} + 3\sq s_{i-1} + 4 s_{i-1} & = & (t_i + s_i\sq) \\
3 t_{i-1} + 4 s_{i-1} + (2 t_{i-1} + 3 s_{i-1}\sq) & = & (t_i + s_i\sq) \\
\end{eqnarray*}
and, by recurrence:
$$ \left\{	\begin{aligned}
t_i & = & 3 t_{i-1} + 4 s_{i-1} \\
s_i & = & 2 t_{i-1} + 3 s_{i-1} 
\end{aligned}	\right.	$$

%\ subsection{Calcolo del limite del rapporto} 
%\ label{limite}
%Osserviamo anche che, al crescere di $i$, il rapporto tra $t$ ed $s$ tende a $\sq$. Questo può essere dimostrato nel seguente modo: se consideriamo la quantità $(3-2\sq)$, identica alla base delle potenze precedenti ma con il segno meno invece che il segno più, e che pertanto possiamo definire \textit{coniugata}, le sue potenze sono ancora le potenze precedenti con il segno cambiato: ovvero, la potenza $i$-esima della coniugata è pari alla coniugata della potenza $i$-esima; la dimostrazione è banale. Tuttavia, a differenza dell'altro caso, il valore $(3-2\sq)$ è minore di $1$, e pertanto $ \lim_{i \rightarrow +\infty} (3-2\sq)^i = \lim_{i \rightarrow +\infty} (t_i - s_i\sq) = 0 $; si potrà quindi scrivere, per ogni $i$, $t_i = k_i s_i$, con $ \lim_{i \rightarrow +\infty} k_i = \sq $. Possiamo anche esprimere la relazione, ponendo $ l_i := k_i - \sq $, come $ t_i = \sq s_i + l_i s_i $, con $ \lim_{i \rightarrow +\infty} l_i = 0 $, in modo tale da avere una quantità che tende a zero. \\
%In questo modo possiamo dimostrare, identificando di che cosa si tratta, che il rapporto tra due numeri sia triangolari che quadrati successivi tende a quella quantità di circa $33.97056$ riscontrata empiricamente. \\

\subsection{Ratio limit: first ratio}
\label{limite}
%Infatti, osserviamo che, poichè $ n = \frac{s}{2} $, allora $ n_i^2 = \frac{s_i^2}{4} $. Vediamo come esprimere $s_i$ unicamente in funzione di $s_{i-1}$ (e non di $t_{i-1}$). Dalle equazioni:
By observing that $ n = \frac{s}{2} $ implies $ n_i^2 = \frac{s_i^2}{4} $, we can express $s_i$ as a function of $s_{i-1}$ and not of $t_{i-1}$. \\
If we define, for every $i$, $t_i = k_i s_i$, $ \lim_{i \rightarrow +\infty} k_i = \sq $ holds (it is straightforward to prove), and we can set $ l_i := k_i - \sq $, so $ t_i = \sq s_i + l_i s_i $, and $ \lim_{i \rightarrow +\infty} l_i = 0 $.
Now, from equations:
$$ \left\{	\begin{aligned}
s_i & = 2 t_{i-1} + 3 s_{i-1} \\
t_i & = \sq s_i + l_i s_i
\end{aligned}	\right.	$$
we obtain: %otteniamo:
\begin{eqnarray*}
s_i & = & 2 \sq s_{i-1} + 2 l_i s_{i-1} + 3 s_{i-1} \\
s_i & = & s_{i-1} (3 + 2\sq + 2 l_i) 
\end{eqnarray*}
from which: %da cui:
\begin{eqnarray*}
s_i^2 & = & s_{i-1}^2 (3 + 2\sq + 2 l_i)^2 \\
4n_i^2 & = & 4n_{i-1}^2 (3 + 2\sq + 2 l_i)^2 \\
n_i^2 & = & n_{i-1}^2 (3 + 2\sq + 2 l_i)^2
\end{eqnarray*}
%Abbiamo detto che, al tendere di $i$ all'infinito, la quantità di $l_i$ tende a $0$. Dalla definizione di limite di una successione, posto un qualsiasi margine di errore come accettabile, esisterà un $i$ sufficientemente grande tale per cui per esso e per tutti i successivi si possa scrivere, approssimando nel margine di errore richiesto:
%$$ n_i^2 = n_{i-1}^2 (3 + 2\sq)^2 $$
%Si verifica che $ (3 + 2\sq)^2 = 17 + 12\sq = (1 + \sq)^4 \cong 33.97056 $, rispettando effettivamente quanto riscontrato nella pratica. \\
%Siamo dunque riusciti ad esprimere quella quantità ricorrente in una forma algebrica.
and, for $ i \rightarrow \infty $:
\begin{eqnarray*}
n_i^2 & = & n_{i-1}^2 (3 + 2\sq)^2 = n_{i-1}^2 (17 + 12\sq) \\ 
      & = & n_{i-1}^2 (1 + \sq)^4 \cong n_{i-1}^2 \cdot 33.97056
\end{eqnarray*}

\subsection{Ratio limit: second ratio, first method}
We prove now in two ways that, if we define:
$$ a_{2,j} = \frac{a_{1,j+1}}{a_{1,j}}-\frac{a_{1,j}}{a_{1,j-1}} $$
then also the ratio $ a_{2,j-1}/a_{2,j} $ tends at the same value for diverging $j$. \\
Here is the first one. \\
We will write alternatively $a_{1,j}$ or $a_j$ for the $j$-th term of the OEIS sequence \texttt{A001110} (see also \cite{asiru,beiler,catarino,dickson,gardner,petkovic,silverman,sloane1,sloane2} and some references therein). \\
We have:
$$	\lim_{j \rightarrow +\infty} \frac{a_{2,j-1}}{a_{2,j}}
=	\lim_{j \rightarrow +\infty} \frac
	{\frac{a_{1,j}}{a_{1,j-1}}-\frac{a_{1,j-1}}{a_{1,j-2}}}
	{\frac{a_{1,j+1}}{a_{1,j}}-\frac{a_{1,j}}{a_{1,j-1}}}
= 	\lim_{j \rightarrow +\infty} \frac
	{\frac{a_{j}}{a_{j-1}}-\frac{a_{j-1}}{a_{j-2}}}
	{\frac{a_{j+1}}{a_{j}}-\frac{a_{j}}{a_{j-1}}} := L_2 $$
Since $ a_j = \frac{s_j^2}{4} $, where $s_j$ is the $j$-th value of $s$ which is solution, for a certain value of $t$ (namely $t_j$), of $t^2-2s^2=1$, we can operate a substitution, implicitely simplifying a $4$ in every fraction:
$$ 	{L_2} = \lim_{j \rightarrow +\infty} \frac
	{\frac{s_{j}^2}{s_{j-1}^2}-\frac{s_{j-1}^2}{s_{j-2}^2}}
	{\frac{s_{j+1}^2}{s_{j}^2}-\frac{s_{j}^2}{s_{j-1}^2}} $$
Now is:
\begin{eqnarray*}
s_{j+1}^2 	& = & s_j^2 \cdot (3 + 2\sq + 2 l_{j+1})^2 		\\
			& = & s_{j-1}^2 \cdot (3 + 2\sq + 2 l_{j+1})^2 
							\cdot (3 + 2\sq + 2 l_j)^2 		\\
			& = & s_{j-2}^2 \cdot (3 + 2\sq + 2 l_{j+1})^2 
							\cdot (3 + 2\sq + 2 l_j)^2
							\cdot (3 + 2\sq + 2 l_{j-1})^2 
\end{eqnarray*}
\begin{eqnarray*}
s_j^2 & = & s_{j-1}^2 \cdot (3 + 2\sq + 2 l_j)^2 			\\
& = & s_{j-2}^2 \cdot (3 + 2\sq + 2 l_j)^2 \cdot (3 + 2\sq + 2 l_{j-1})^2
\end{eqnarray*}
$$ s_{j-1}^2 = s_{j-2}^2 \cdot (3 + 2\sq + 2 l_{j-1})^2 $$
where $ l_j = t_j/s_j - \sq $, and $ l_j \rightarrow 0 $ for $ j \rightarrow +\infty $. \\
This lead to the ratios:
\begin{eqnarray*}
\frac{s_{j+1}^2}{s_j^2} & = & 
%\frac{s_{j-2}^2 
%\cdot(3+2\sq+2l_{j+1})^2 \cdot(3+2\sq+2l_j)^2 \cdot(3+2\sq+2 l_{j-1})^2}
%{s_{j-2}^2 \cdot (3 + 2\sq + 2 l_j)^2 \cdot (3 + 2\sq + 2 l_{j-1})^2}
(3 + 2 \sq +2 l_{j+1})^2 \\
\frac{s_j^2}{s_{j-1}^2} & = &
%\frac{s_{j-2}^2 \cdot (3 + 2\sq + 2 l_j)^2 \cdot (3 + 2\sq + 2 l_{j-1})^2}
%{s_{j-2}^2 \cdot (3 + 2\sq + 2 l_{j-1})^2}
(3 + 2 \sq +2 l_j)^2 \\
\frac{s_{j-1}^2}{s_{j-2}^2} & = & (3 + 2 \sq +2 l_{j-1})^2 \\
\end{eqnarray*}
We can now rewrite $ L_2 $ by using the ratios:
$$ 	L_2 = \lim_{j \rightarrow +\infty} 
	\frac	{ (3 + 2\sq + 2l_j)^2 - (3 + 2\sq + 2l_{j-1})^2 }
			{ (3 + 2\sq + 2l_{j+1})^2 - (3 + 2\sq + 2l_j)^2 } $$
Now the square differences can be rewritten as a product of a sum and a difference:
$$ 	L_2 = \lim_{j \rightarrow +\infty}
	\frac	{ (6 + 4\sq + 2(l_{j-1} + l_j)) \cdot 2(l_{j-1} - l_j)}
			{ (6 + 4\sq + 2(l_j + l_{j+1})) \cdot 2(l_j - l_{j+1})} $$
Considering the fact that $l_j$ tends to zero for diverging $j$, we can both approximate $ 6 + 4\sq + 2(l_{j-1} + l_j) $ and $ 6 + 4\sq + 2(l_j + l_{j+1}) $ with $ 6 + 4\sq $. Then: \\
$$ 	L_2 = \lim_{j \rightarrow +\infty}
\frac{(6 + 4\sq) \cdot 2(l_{j-1} - l_j)}{(6 + 4\sq) \cdot 2(l_j - l_{j+1})}
= \lim_{j \rightarrow +\infty} \frac{l_{j-1} - l_j}{l_j - l_{j+1}} $$ 
and so:
$$ L_2 = \lim_{j \rightarrow +\infty} \frac{l_{j-1} - l_j}{l_j - l_{j+1}}
	   = \lim_{j \rightarrow +\infty} 
	     \frac{(k_{j-1}-\sq) - (k_j-\sq)}{(k_j-\sq) - (k_{j+1}-\sq)} 
	   = \lim_{j \rightarrow +\infty} \frac{k_{j-1} - k_j}{k_j - k_{j+1}} $$
where $ k_j = t_j/s_j $.
$$ L_2 = \lim_{j \rightarrow +\infty}
		 \frac	{\frac{t_{j-1}}{s_{j-1}} - \frac{t_j}{s_j}}
		 		{\frac{t_j}{s_j} - \frac{t_{j+1}}{s_{j+1}}} 
	   = \lim_{j \rightarrow +\infty}
	   	 \frac	{\frac{t_{j-1} s_j - t_j s_{j-1}}{s_j s_{j-1}}}
	   	 		{\frac{t_j s_{j+1} - t_{j+1} s_j}{s_j s_{j+1}}}
	   = \lim_{j \rightarrow +\infty}
	     \frac	{(t_{j-1} s_j - t_j s_{j-1}) \cdot s_j \cdot s_{j+1}}
	     		{(t_j s_{j+1} - t_{j+1} s_j) \cdot s_j \cdot s_{j-1}}
$$
By proceeding with calculations we can state:
$$ L_2 = \lim_{j \rightarrow +\infty} \left(
		 \frac{s_{j+1}}{s_{j-1}} \cdot
	   	 \frac{t_{j-1} s_j - t_j s_{j-1}}{t_j s_{j+1} - t_{j+1} s_j} \right)
	   = (3+2\sq)^2 \cdot \lim_{j \rightarrow +\infty}
	   	 \frac{t_{j-1} s_j - t_j s_{j-1}}{t_j s_{j+1} - t_{j+1} s_j}
$$
where $ s_{j+1}/s_{j-1} = (s_{j+1}/s_j) \cdot (s_j/s_{j-1}) $, and the limit of both factors is equal to $(3+2\sq)$. \\
For the remaining limit, we consider just the denominator:
\begin{eqnarray*}
		t_j s_{j+1} - t_{j+1} s_j 
	& = &   t_j (2 t_j + 3 s_j) - (3 t_j + 4 s_j) s_j
    = 	2 t_j^2 + 3 s_j t_j - 3 s_j t_j - 4 s_j^2  \\
    & = & 	2 t_j^2 - 4 s_j^2 
 	= 	2 t_j^2 - 4 s_j^2 = 2 (t_j^2 - 2 s_j^2) = 2 \cdot 1 = 2
\end{eqnarray*}
where the factor in brackets is equal to $1$ for every $j$, because $(t_j,s_j)$ is a solution of the Pell equation $ t_j^2 - 2 s_j^2 = 1 $. In particular the same result is obtaining by considering the numerator, because it is just the denominator with indices shifted by one. Then the ratio is constant and equal to $1$; so is the limit for $j \rightarrow 0 $, and:
$$ L_2 = (3+2\sq)^2 = (17+12\sq) = (1+\sq)^4 $$
as we wanted to prove.

\subsection{Ratio limit: second ratio, second method}
We will see now an alternate way to get that result. \\
We know the solutions of the Pell equation to be $ t_j + s_j \sq = (3+2\sq)^j $, and also that $ t_j - s_j \sq = (3-2\sq)^j $. Observed $ (3-2\sq)=(3+2\sq)^{-1} $, and defined $ \beta:=(1+\sq)$, hence $ (3+2\sq) = \beta^2 $, $ (3-2\sq) = \beta^{-2} $, by respectively summing and subtracting:
$$ \left\{	\begin{aligned}
t_j = \frac{\beta^{2j} + \beta^{-2j}}{2} \\
s_j = \frac{\beta^{2j} - \beta^{-2j}}{2\sq}
\end{aligned}	\right.	$$
%\begin{eqnarray*}
%t_j = \frac{\beta^{2j} + \beta^{-2j}}{2} \\
%s_j = \frac{\beta^{2j} - \beta^{-2j}}{2\sq}
%\end{eqnarray*}
This allows us to write a closed formula, from which we can generate numbers which are both triangulars and squares:
$$ a_{1,j} 	= \frac{s_j^2}{4} 	= \frac{\beta^{4j} + \beta^{-4j} - 2}{32} 
								= \frac{\alpha^{j} + \alpha^{-j} - 2}{32} $$
by setting $ \alpha = \beta^4 = (1+\sq)^4 $. \\
We can obtain via these calculations the well-known result:
$$ 	\lim_{j \rightarrow +\infty} \frac{a_{1,j}}{a_{1,j-1}}
	= 	\lim_{j \rightarrow +\infty}
		\frac	{\frac{\alpha^{j} + \alpha^{-j} - 2}{32}}
				{\frac{\alpha^{j-1} + \alpha^{1-j} - 2}{32}} 
	=	\lim_{j \rightarrow +\infty}
		\frac{\alpha^{j}+\alpha^{-j}-2}{\alpha^{j-1}+\alpha^{1-j}-2} 
	=	\lim_{j \rightarrow +\infty} \frac{\alpha^j}{\alpha^{j-1}}
	=	\alpha $$
considering that $ |\alpha| > 1 $ and so other terms are trascurable for $ j \rightarrow +\infty $. \\
In an analogue way we can compute:
$$  \lim_{j \rightarrow +\infty} \frac{a_{2,j-1}}{a_{2,j}}
=	\lim_{j \rightarrow +\infty} \frac
	{\frac{a_{1,j}}{a_{1,j-1}}-\frac{a_{1,j-1}}{a_{1,j-2}}}
	{\frac{a_{1,j+1}}{a_{1,j}}-\frac{a_{1,j}}{a_{1,j-1}}} 
=	\lim_{j \rightarrow +\infty} \frac
	{\frac{\alpha^{j}+\alpha^{-j}-2}{\alpha^{j-1}+\alpha^{1-j}-2}
	-\frac{\alpha^{j-1}+\alpha^{1-j}-2}{\alpha^{j-2}+\alpha^{2-j}-2}}
	{\frac{\alpha^{j+1}+\alpha^{-1-j}-2}{\alpha^{j}+\alpha^{j}-2}
	-\frac{\alpha^{j}+\alpha^{j}-2}{\alpha^{j-1}+\alpha^{1-j}-2}} $$
by implicitely simplifying a $32$ in every fraction. \\
The use of standard algebra techniques gives the subsequent results.
$$  \lim_{j \rightarrow +\infty} \frac
	{\frac{\alpha^{j}+\alpha^{-j}-2}{\alpha^{j-1}+\alpha^{1-j}-2}
	-\frac{\alpha^{j-1}+\alpha^{1-j}-2}{\alpha^{j-2}+\alpha^{2-j}-2}}
	{\frac{\alpha^{j+1}+\alpha^{-1-j}-2}{\alpha^{j}+\alpha^{j}-2}
	-\frac{\alpha^{j}+\alpha^{j}-2}{\alpha^{j-1}+\alpha^{1-j}-2}} 
=	\lim_{j \rightarrow +\infty} \frac
	{\frac
	{(\alpha^{j}+\alpha^{-j}-2) \cdot (\alpha^{j-2}+\alpha^{2-j}-2) 
	 \minusm (\alpha^{j-1}+\alpha^{1-j}-2)^2}
	{(\alpha^{j-1}+\alpha^{1-j}-2) \cdotm (\alpha^{j-2}+\alpha^{2-j}-2)}}
	{\frac
	{(\alpha^{j+1}+\alpha^{-1-j}-2) \cdot (\alpha^{j-1}+\alpha^{1-j}-2) 
	 \minusm (\alpha^{j}+\alpha^{-j}-2)^2}
	{(\alpha^{j}+\alpha^{-j}-2) \cdotm (\alpha^{j-1}+\alpha^{1-j}-2)}} $$
By rearranging:
{\footnotesize
$$  \lim_{j \rightarrow +\infty} 
	\left( \frac
	{(\alpha^{j}+\alpha^{-j}-2) \cdot (\alpha^{j-2}+\alpha^{2-j}-2) 
	 \minusm (\alpha^{j-1}+\alpha^{1-j}-2)^2}
	{(\alpha^{j+1}+\alpha^{-1-j}-2) \cdot (\alpha^{j-1}+\alpha^{1-j}-2) 
	 \minusm (\alpha^{j}+\alpha^{-j}-2)^2}
	 \cdot \frac
	{(\alpha^{j}+\alpha^{-j}-2) \cdotm (\alpha^{j-1}+\alpha^{1-j}-2)}
	{(\alpha^{j-1}+\alpha^{1-j}-2) \cdotm (\alpha^{j-2}+\alpha^{2-j}-2)}
	\right) $$
}
by semplifying in the right factor:
{\footnotesize
$$  \lim_{j \rightarrow +\infty} 
	\left( \frac
	{(\alpha^{j}+\alpha^{-j}-2) \cdot (\alpha^{j-2}+\alpha^{2-j}-2) 
	 \minusm (\alpha^{j-1}+\alpha^{1-j}-2)^2}
	{(\alpha^{j+1}+\alpha^{-1-j}-2) \cdot (\alpha^{j-1}+\alpha^{1-j}-2) 
	 \minusm (\alpha^{j}+\alpha^{-j}-2)^2}
	 \cdot \frac
	{(\alpha^{j}+\alpha^{-j}-2)}{(\alpha^{j-2}+\alpha^{2-j}-2)}
	\right) $$
}
and by explicitely calculating the left factor:
{\footnotesize
$$  \lim_{j \rightarrow +\infty} \left( \frac
	{ \alpha^2 - 2\alpha^j + \alpha^{-2} - 2\alpha^{-j} - 2\alpha^{j-2}
	+2\alpha^{2-j} - 2 + 4\alpha^{j-1} + 4\alpha^{1-j} }
	{ \alpha^2 - 2\alpha^{j+1} + \alpha^{-2} - 2\alpha^{-1-j}
	- 2\alpha^{j-1}	+ 2\alpha^{1-j} - 2 + 4\alpha^{j} + 4\alpha^{-j} }
	\cdot \frac
	{(\alpha^{j}+\alpha^{-j}-2)}{(\alpha^{j-2}+\alpha^{2-j}-2)} 
	\right) $$	
}
Considered the fact that there is the limit operator, we can consider just the elements depending on $j$, in which the coefficents of them at the exponential are positive, because the others are trascurabile with respect to them, considering the operations we are doing. This finally gives:
$$ 	\lim_{j \rightarrow +\infty} \left( \frac
	{ 	-2\alpha^j - 2\alpha^{j-2} + 4\alpha^{j-1}		}
	{	-2\alpha^{j+1} - 2\alpha^{j-1} + 4\alpha^{j}	}  
	\cdot \frac {\alpha^j}{\alpha^{j-2}} \right) 
=	(\alpha^{-1} \cdot \alpha^2) = \alpha $$
again, as we wanted to prove.

\section{Open points}
We conjecture that the result holds for every $h$-th ratio, $ h \geq 3 $, defined by:
$$ a_{h,j} = \frac{a_{h-1,j+1}}{a_{h-1,j}}-\frac{a_{h-1,j}}{a_{h-1,j-1}} $$
This means that it holds:
$$ \lim_{j \rightarrow +\infty} \frac{a_{h,{j-1}}}{a_{h,j}} 
=  \lim_{j \rightarrow +\infty} \frac
   {\frac{a_{h-1,j}}{a_{h-1,j-1}}-\frac{a_{h-1,j-1}}{a_{h-1,j-2}}}
   {\frac{a_{h-1,j+1}}{a_{h-1,j}}-\frac{a_{h-1,j}}{a_{h-1,j-1}}}
= (1+\sq)^4 $$
but we are not able to either prove or disprove it, at the moment. \\
On the other hand, it can be investigated whether similar results can be written for other sequences of integers figurate in more than one way, like both triangular and pentagonal, both square and pentagonal, and so on.

\end{document}